\documentclass[12pt]{article}
\usepackage{amssymb, amsmath}
\newtheorem{thm}{Theorem}  
\newtheorem{defin}{Definition} 

\newcommand{\n}{\nonumber}
\newcommand{\bb}{\begin{equation}}
\newcommand{\ee}{\end{equation}}
\newcommand{\bq}{\begin{eqnarray}}
\newcommand{\eq}{\end{eqnarray}}
\newcommand{\bqn}{\begin{eqnarray*}}
\newcommand{\eqn}{\end{eqnarray*}}

\renewcommand{\a}{\alpha}

\renewcommand{\div}{\mathrm{div}\, }

\newcommand{\s}{\sigma}
\newcommand{\si}{\sigma_R (x)}

\newcommand{\vare}{\varepsilon}

\newcommand{\bu}{\mbox{\boldmath $u$}}

\newcommand{\bm}{\mbox{\boldmath $m$}}
\newcommand{\bx}{\mbox{\boldmath $x$}}
\newcommand{\bpsi}{\mbox{\boldmath $\psi$}}
\newcommand{\bphi}{\mbox{\boldmath $\phi$}}

\begin{document}
\title{Blow-up, zero $\alpha$ limit and the Liouville type theorem
for   the Euler-Poincar\'{e} equations}
\author{Dongho Chae$^{*}$ and Jian-Guo Liu$^{\dagger}$}
 \date{}
\maketitle

\begin{center}
${}^{(*)}$Department of Mathematics\\
  Sungkyunkwan University\\
 Suwon 440-746, Korea\\
 email: chae@skku.edu
\end{center}

\begin{center}
${}^{(\dagger)}$Department of Physics and Department of Mathematics \\
Duke University \\
Durham, NC 27708, USA \\
email: jliu@phy.duke.edu
\end{center}

\begin{abstract}
In this paper we study the Euler-Poincar\'{e} equations in $\Bbb
R^N$. We prove local existence of weak solutions in $W^{2,p}(\Bbb
R^N),$ $p>N$, and local existence of unique classical solutions in
$H^k (\Bbb R^N)$, $k>N/2+3$, as well as a blow-up criterion. For the zero dispersion equation($\alpha=0$) we prove a finite time
blow-up of the classical solution. We
also prove that as the dispersion parameter vanishes, the weak solution converges to a solution of the zero dispersion equation with sharp rate as $\alpha\to0$,  provided that the limiting solution belongs to $C([0, T);H^k(\Bbb R^N))$ with
$k>N/2 +3$.   For the {\em stationary weak
solutions} of the Euler-Poincar\'{e} equations we prove a Liouville
type theorem. Namely, for $\alpha>0$ any weak solution $\mathbf{u}\in  H^1(\Bbb
R^N)$ is $\mathbf{u}=0$; for $\alpha=0$ any weak solution $\mathbf{u}\in  L^2(\Bbb
R^N)$ is $\mathbf{u}=0$.
\end{abstract}


\vskip 0.4cm

\section{Introduction}
We consider the following Euler-Poincar\'{e} equations in $\Bbb R^N$:
$$
(EP) \left\{ \aligned & \partial_t \bm +(\bu\cdot \nabla ) \bm
+(\nabla \bu)^\top \bm
+(\mathrm{div}\, \bu)\bm=0,\\
 &\bm=(1-\alpha \Delta )\bu ,\\
 &\bu_0 (x)=\bu_0 \,,
 \endaligned \right.
 $$
 where $\bu: \Bbb R^N \to \Bbb R^N$ is the velocity, $\bm: \Bbb R^N \to \Bbb R^N$ represents the momentum,
 constant $\sqrt{\alpha}$ is a length scale parameter, $(\nabla \bu)^\top=$ the transpose of $(\nabla \bu)$.
 The Euler-Poincar\'{e} equations  arise in diverse scientific applications and enjoy several remarkable
properties both in the one-dimensional and multi-dimensional cases.

 The  Euler-Poincar\'{e} equations were first studied 
 by Holm, Marsden, and Ratiu  in 1998 as  a framework for modeling and analyzing fluid dynamics \cite{HMR1, HMR2}, particularly
 for nonlinear shallow water waves,  geophysical fluids and turbulence modeling.
There are intensive researches on analogs viscous or inviscid, incompressible Lagrangian averaged models. We refer  to \cite{ Chen, Foias, Marsden} for results on  Navier-Stokes-$\alpha$ model in terms of existence and uniqueness,  zero $\alpha$ limit to the Navier-Stokes equations, global attractor, etc.
We refer to \cite{Bardos, HNP, TitiXin} for results on analysis and simulation of vortex sheets with Birkhoff-Rott-$\alpha$ or  Euler-$\alpha$ approximation.

For one-dimension, the Euler-Poincar\'{e} equations coincide with the dispersion-less case of  Camassa-Holm (CH) equation \cite{CH}:
$$
(CH) \quad \partial_t m + 2 u \partial_x m + \partial_x u \, m = 0, \qquad m=(1 - \alpha \partial_{xx} ) u \,.
$$
The solutions  to (CH) are characterized by a discontinuity in the first derivative at their peaks
and are thus referred to as peakon solutions. (CH) is completely integrable with a bi-Hamiltonian structure and their peakon
solutions are true solitary waves that emerge from the initial data. Peakons exhibit a remarkable stability--their
identity is preserved through nonlinear interactions, see, e.g. \cite{CH,  HSS}.
There are many comprehensive analysis on (CH)  in the literature.
 We refer to
a review paper \cite{MoLu2004} for a survey of recent results on well-poseness and existence of local and global weak solutions for
(CH).
The existence of a global weak solution and uniqueness was proven in
\cite{BrCo2007,  Chertock, CM, CE98a,XZ01}. A class of the so called weak-weak solution was studied in \cite{XZ01}.
The breakdown of the solution for (CH) was studied in \cite{McKean}.

The Euler-Poincar\'{e} equations have many further interpretations beyond fluid applications. For instance,
in 2-D, it is exactly the same as  the averaged template matching equation for computer vision (see, e.g.,
\cite{HiMaAr2001,HoMa2005,HoRaTrYo2004}). The Euler-Poincar\'{e} equations  also has  important applications in computational anatomy (see, e.g, \cite{HSS, Younes}).
 The  Euler-Poincar\'{e} equations can also be regarded as an evolutionary equation for a geodesic motion on a diffeomorphism group  and it is associated with  Euler-Poincar\'{e}  reduction via symmetry \cite{Arnold, EbinMarsden, HSS, Khesin,Younes}.
  We refer to a recent book \cite{HSS} for a comprehensive review on the subject.

  The organization of the paper is as follows.
  In Section 2, we give some preliminary discussions of the Euler-Poincar\'{e}  equations and we state a theorem on  local existence of weak solution in $W^{2,p}(\Bbb
R^N),$ $p>N$, and local existence of unique classical solutions in
$H^k (\Bbb R^N)$, $k>N/2+3$.

In Section 3, we prove a theorem on a  blow-up criterion, as well as,  a theorem on finite time
blow-up of the classical solution for the zero dispersion equation. For classic solutions with reflection symmetry, the divergence $\nabla \cdot \bu$ satisfy a Riccati  equation at the invariant point under the reflection transformation  and hence there is a finite time blow up if  the divergence is initially negative. 

 In Section 4, we
 prove that as the dispersion parameter $\alpha$ vanishes, the weak solution converges to a solution of the zero dispersion equation with a sharp rate as $\alpha \to 0$,  provided that the limiting solution belongs to $C([0, T);H^k(\Bbb R^N))$ with
$k>N/2 +3$.

Finally, for the {\em stationary weak
solutions} of the Euler-Poincar\'{e} equations we prove a Liouville
type theorem in Section 5. For $\alpha>0$, we prove that any weak solution $\mathbf{u}\in  H^1(\Bbb
R^N)$ is  $\mathbf{u}=0$. For $\alpha=0$, any weak solution $\mathbf{u}\in  L^2(\Bbb
R^N)$ is $\mathbf{u}=0$.
This is a surprising result, 
as all the previous  Liouville
type results are for dissipative systems. This is the first  Liouville
type theorem for non-dissipative systems.

  We also give a proof of the local existence and uniqueness theorem in  Appendix.

\section{Preliminaries and local existence}

In this section, we discuss some mathematical structures of (EP)  and then we state a local existence theorem for the weak solution and the classic solution.
We refer  to \cite{HoMa2005, HSS} for more in-depth discussions on (EP).

(EP) can be recast as
\bb
 \partial_t \bm + \nabla\cdot (\bu\otimes \bm)+(\nabla \bu)^\top \bm =0 \,.
\ee
 The last term  above can  be written in a conservative/tensor form
 \bqn
\sum_{j=1}^N \partial_i u_j m_j&=& \sum_{j=1}^N \partial_iu_j u_j
-\a \sum_{j,k=1}^N \partial_i u_j
\partial_k ^2
u_j\\
&=& \frac12  \partial_i |\bu|^2 -\a \sum_{j,k=1}^N \partial_k
(\partial_i u_j
\partial_k u_j )+\a \sum_{j,k=1}^N \partial_k \partial_i u_j \partial_k u_j \\
&=&\frac12  \partial_i |\bu|^2-\a \sum_{j,k=1}^N
\partial_j(\partial_i
u_k \partial_j u_k ) +\frac{\a}{2} \sum_{j,k=1}^N \partial_i (\partial_k u_j)^2\\
&=& \sum_{j=1}^N \partial_j\left( \frac12\delta_{ij} |\bu|^2 -\a
\partial_i \bu\cdot
\partial_j \bu +\frac{\a}{2}\delta_{ij} |\nabla \bu|^2 \right).
\eqn
Set stress-tensor
 $$
 T_{ij} = m_i u_j +\frac{\delta_{ij} }{2} |\bu|^2 -\a \partial_i \bu\cdot
\partial_j \bu +\frac{\a\delta_{ij}}{2} |\nabla \bu|^2 .
$$
Then (EP) becomes
 \bb \label{tensor}
 \partial_t m_i + \sum_{j=1}^N \partial_j T_{ij} = 0 \,.
 \ee
The first term in $T_{ij} $ involves a second order derivative of $\bu$ and it can be rewritten as
$$
 m_i u_j = u_i u_j + \alpha\sum_{k=1}^{N} \partial_k u_{i} \partial_k u_{j} -  \alpha \sum_{k=1}^{N} \partial_k \big(u_{j} \partial_k u_{i} \big).
$$
The symmetric part of tensor $T$ is given by
 \begin{equation}\label{Ta}
 T^a= \bu \otimes \bu  + \a \nabla \bu \nabla \bu^\top
  -\a \nabla \bu^\top \nabla \bu
+ \frac12 (|\bu|^2+\a  |\nabla \bu|^2) \mbox{ Id }
\end{equation}
and the remainder terms in $T$ are given by
\bb \label{Tb}
T^{b}_{i,j} = - \alpha \sum_{k=1}^{N} \partial_k \big( u_{j} \partial_k u_{i} \big) \,.
\ee
Hence $T =  T^a + T^b$.
In view of this, the natural definition of the weak solution of (EP) would be:
\begin{defin} \quad
$\bu\in 
L^\infty(0, T;
H^{1}_{loc} (\Bbb R^N))$ is a weak solution of (EP) with initial
data $\bu_0 \in H^{1}_{loc} (\Bbb R^N)$ if the following equation
holds for all vector field $\bphi (x,t)$ such that $\bphi(\cdot ,t)\in C_0
^\infty (\Bbb R^N) $ for all $ t\in [0, T)$ and $\bphi (x, \cdot)\in
C^1_0([0, T))$ for all $ x\in \Bbb R^N$
 \bq \label{weak}
 \lefteqn{ \int_0 ^T \int_{\Bbb R^N}
  \bigl(  \bu \cdot \bphi_t +   \alpha \nabla \bu : \nabla \bphi_t  \bigr) dxdt
 +\int_{\Bbb R^N}
\bigl(  \bu_0 \cdot \bphi(\cdot,0) +  \alpha \nabla \bu_0 : \nabla \bphi(\cdot,0) \bigr)     \,dx
} \n \\
 &&  +  \int_0 ^T \int_{\Bbb R^N}
 T^{a}:\nabla \bphi (x,t)\, dxdt
 +\a  \sum_{i,j,k=1}^{N} \int_0 ^T \int_{\Bbb R^N} u_j \partial_k u_i \partial_{j} \partial_{k} \phi_i \, dxdt=0,
 \eq
 where $T^{a}$ is given by (\ref{Ta}).
\end{defin}

(EP) also has a natural  Hamiltonian structure. Set
$$
    {\cal H} = \frac12 \int_{\Bbb R^{N}}  \bu\cdot \bm \,d\bx.
$$
then $\frac{ \delta {\cal H} }{\delta \bm }=\bu$ and (EP) can be recast as
\bb \label{Hamiltonian}
   \partial_t \bm = - {\cal A}\frac{ \delta {\cal H} }{\delta \bm },
\ee
where ${\cal A}$ is an anti-symmetric operator defined by
$$
  {\cal A} \bu = \sum_{j=1}^N \partial_j (m_i u_j) + \sum_{j=1}^N \partial_i u_j m_j.
$$
Consequently, from (\ref{tensor}) and  (\ref{Hamiltonian}), there are two conservation laws
$$
  \frac{d}{dt} \int_{\Bbb R^{N}}  \bm \,d\bx = 0, \quad
   \frac{d}{dt}  \int_{\Bbb R^{N}} ( |\bu|^{2} + \a |\nabla \bu|^2 ) \,d\bx = 0.
$$
For the one-dimensional case,  (EP) coincides with the dispersion-less case of  Camassa-Holm (CH) equation
and there is an additional Hamiltonian structure and a Lax-pair which leads to a complete integrability of (CH) \cite{CH}. We refer to \cite{Fokas} for a general discussion on bi-Hamiltonian system and complete integrability.

When $\alpha=0$, the above Hamiltonian structure shows that (EP) is a symmetric hyperbolic system of conservation laws
\begin{equation}\label{hyp}
 \left\{\aligned &\partial_t \bu +  \mathrm{div} (\bu\otimes \bu) + \frac12 \nabla |\bu|^{2} = 0\\
 &\bu(x,0)=\bu_0
 \endaligned\right.
\end{equation}
which possess a global convex entropy function
\begin{equation}\label{entropy}
\frac12 \partial_t |\bu|^2 +  \mathrm{div} (|\bu|^{2}  \bu) = 0 \,.
\end{equation}
We refer (\ref{hyp}) as the zero dispersion equation.
Indeed, we can recast it in a usual form of  a symmetric  hyperbolic system (we state it in $\Bbb R^3$):
$$
   \bu_t + A \bu_x + B \bu_y + C \bu_z = 0
$$
with
$$
  \bu= \begin{pmatrix} u \\ v \\ w \end{pmatrix}, A=\begin{pmatrix} 3u & v & w \\ v & u &0 \\ w &0 &u \end{pmatrix}, \cdots
$$
$A$ is a symmetric matrix and has  three eigenvalues: $ u$, $2u + |\bu|$, $2u - |\bu|$, corresponding to one linearly degenerate field, and two genuinely nonlinear fields, respectively, when $\bu\ne 0$.

We shall remark that although the high dimensional Burgers equation has a similar structure as (\ref{hyp}), it does not   possess  a global convex entropy. In section 5, we will prove a Liouville type theorem for (\ref{hyp}).
This theorem does not hold true for the high dimensional Burgers equation.

Now we introduce some notations and then we state a theorem on local existence of the weak solution and local existence and  uniqueness of the classical solution.

 For $s\in \Bbb R$ and  $p\in [1, \infty]$ we
define the Bessel potential space $L^{s,p} (\Bbb R^N)$ as follows
$$
L^{s,p} (\Bbb R^N)=\{ f\in L^p (\Bbb R^N)\, |\, \|(1-\Delta
)^{\frac{s}{2}} f\|_{L^p}:=\|f\|_{L^{s,p}} <\infty\}.
$$
For $s\in \Bbb N\cup\{0\}$ it is well-known that $L^{s,p} (\Bbb
R^N)$ is equivalent to the standard Sobolev space $W^{s,p} (\Bbb
R^N)$(see e.g. \cite{ste}). This, in turn, implies immediately that
there exist $C_1, C_2$ such that
 \bb
C_1 \|\bu\|_{W^{k+2,p}} \leq \|\bm\|_{L^{k,p}} \leq C_2
\|\bu\|_{W^{k+2,p}}
 \ee
for all $k\in \Bbb N\cup \{ 0\}, p\in (1, \infty)$.
As usual we denote $H^{s}(\Bbb R^N)=W^{s, 2}(\Bbb R^N)$.

\begin{thm} \label{thm1}
\begin{itemize}
\item[(i)] Assume $\alpha >0$ and $\bu_0 \in W^{2,p} (\Bbb R^N)$ with $p>N$. Then, there
exists $T=T(\|\bu_0 \|_{W^{2,p}})$ such that a weak solution to (EP)
exists, and belongs to $\bu \in L^\infty(0, T;W^{2,p} (\Bbb R^N))\cap Lip(0, T; W^{1,p} (\Bbb R^N))$.
\item[(ii)]
Let $\alpha >0$ and  $\bu_0 \in H^k (\Bbb R^N)$ with $k>N/2+3$. Then, there exists
$T=T(\|\bu_0 \|_{H^k})$ such that a classic solution to (EP) exists
uniquely, and belongs to $\bu \in C([0, T); H^k (\Bbb R^N))$.
\item[(iii)]
For $\alpha=0$,  (EP) is a symmetric hyperbolic system of conservation laws with a  convex entropy. Consequently, if
  $\bu_0 \in H^k (\Bbb R^N)$ with $k>N/2+1$. Then, there exists
$T=T(\|\bu_0 \|_{H^k})$ such that a classic solution to (EP) exists
uniquely, and belongs to $\bu \in C([0, T); H^k (\Bbb R^N))$.
  \end{itemize}
\end{thm}
The proof of  symmetric hyperbolicity  and existence of convex entropy in (iii) are given in (\ref{hyp})-(\ref{entropy}).
The proof of existence of the unique classic solution for symmetric hyperbolic system is  standard, see e.g \cite{Majda}.
The proof of (i) and (ii) is also rather standard and will be given in the Appendix for completeness.

\section{Finite time blow up}

 In this section, we first present a theorem on a blow-up criterion and then we prove a theorem on finite time blow up for the zero dispersion equation.

 We denote the deformation tensor
 for $\bu$ by $S=(S_{ij})$, where $S_{ij}:= \frac12 (\partial_i u_j +\partial_j u_i )$.
  We recall the Besov space $\dot{B}^0_{\infty,
\infty}$, which is defined as follows. Let $\{\psi_m\}_{m\in \Bbb
Z}$ be the Littlewood-Paley partition of unity, where the Fourier
transform $\hat{\psi}_m (\xi)$ is supported on the annulus $\{ \xi
\in \Bbb R^N\, |\, 2^{m-1}\leq |\xi|< 2^{m}\}$(see e.g. \cite{tay}).
Then,
$$
f\in \dot{B}^0_{\infty, \infty}\quad \mbox{if and
  only if }\quad \sup_{m\in \Bbb Z} \|\psi_m *f\|_{L^\infty} :=
  \|f\|_{\dot{B}^0_{\infty,
\infty}}<\infty.
$$
The following is a well-known embedding result,
\bb
L^\infty(\Bbb
R^N)\hookrightarrow BMO(\Bbb R^N) \hookrightarrow\dot{B}^0_{\infty,
\infty} (\Bbb R^N).
\ee

\begin{thm} \label{thm2}
 For $\alpha \ge 0$, we have the following finite time blow-up criterion of the
local solution of (EP) in $\bu \in C([0, t_*); H^k (\Bbb R^N))$,
$k>N/2+3$.
  \bb \label{BlowUpCriterion}
  \lim\sup_{t\to t_*} \|\bu(t)\|_{H^k} =\infty \quad \mbox{if and
  only if }\quad \int_0 ^{t_*} \|S (t)\|_{\dot{B}^0_{\infty, \infty}} dt =\infty \,.
  \ee
\end{thm}
{\it Remark 1.1 } Combining the embedding relation, $W^{1,N} (\Bbb
R^N)\hookrightarrow BMO (\Bbb R^N)\hookrightarrow\dot{B}^0_{\infty,
\infty} (\Bbb R^N)$ with the inequality $\|D^2 \bu\|_{L^p} \leq
C\|\bm \|_{L^p} $ for $p\in (1, \infty)$(see (\ref{pseudo}) below), we have
$$\|S \|_{\dot{B}^0_{\infty, \infty}} \leq C \|S\|_{BMO}\leq C \|D S\|_{L^N} \leq C \|D^2 \bu\|_{L^N} \leq C \|\bm
\|_{L^N}.$$
 Therefore we obtain the following criterion as an
immediate corollary of the above theorem: for all $p>N$
 \bb \lim\sup_{t\to t_*}
\|\bm(t)\|_{L^p} =\infty \quad \mbox{if and
  only if }\quad \int_0 ^{t_*} \|\bm (t)\|_{L^N} dt =\infty.
\ee
\noindent{\it Remark 1.2 } In the one dimensional case of
the Camassa-Holm equation (CH) the above criterion implies that  finite
time blow-up does not happen  if $\int_0 ^{t} \|\bu_{xx}
(\tau)\|_{L^1}d\tau <\infty $ for all $t>0$. Thanks to the
conservation law we have $\sup_{0<\tau <t} \|\bu_x (\tau)\|_{L^2}
\leq \|\bu_0 \|_{H^1}< \infty$ for all $t>0$. Since we have
embedding $W^{2,1} (\Bbb R) \hookrightarrow H^{1} (\Bbb R)$, and
we do have finite time blow-up for (CH) \cite{McKean}, our
criterion is sharp in
this one dimensional case.
\\
\ \\
\noindent{\it Proof of Theorem 2  }
We only give a proof for the case $\alpha>0$.
The proof for the case $\alpha=0$ is similar and simpler hence will be omitted.

Using estimates (\ref{I123}, \ref{J1}, \ref{J2}, \ref{I23}) for $I_1, I_2, I_3$ in the proof of Theorem 1 in the Appendix,
one has
 \bqn
 \frac{d}{dt} \|\bm (t) \|_{H^k}&\leq&  C(\|\nabla \bu\|_{L^\infty} + \|\bm
\|_{L^\infty} +\|\nabla \bm \|_{L^\infty})\|\bm (t) \|_{H^k}\\
 &\leq & C(\|\bm \|_{L^p} + \|D\bm \|_{L^p} +
 \|D^2 \bm \|_{L^p})\|\bm (t) \|_{H^k} \,.
 \eqn
Hence,
\bb\label{BUC_Hk}
 \|\bm (t) \|_{H^k}\leq \|\bm _0\|_{H^k} \exp\left[ C\int_0 ^t
 \left\{\|\bm (\tau)\|_{L^p} + \|D\bm (\tau)\|_{L^p}
  + \|D^2 \bm (\tau)\|_{L^p}\right\}d\tau \right]
\ee
 for $k>N/2+1$ and $p>N$, where we used the Sobolev embedding.
Consequently,  blow up of  $\|\bm (t) \|_{H^k}$ as $t\to t^{*}$ implies that at least one of $\|\bm (t)\|_{L^p}$, $\|D\bm (t)\|_{L^p}$
 and $\|D^2 \bm (t)\|_{L^p}$ blow up as $t\to t^{*}$. In the following three steps, we show  blow-up criterion for each of them are all given by  (\ref{BlowUpCriterion}).

{\it Step 1}. \quad
We first recall the following logarithmic Sobolev inequality(see
e.g. \cite{tay}),
 \bb\label{log}
  \|f\|_{L^\infty} \leq C(1+\|f\|_{\dot{B}^0_{\infty, \infty}}
  )(\log (1+ \|f\|_{W^{s,p}} )),
  \ee
  where $s>0, 1<p<\infty$ and $ sp>N$. From the estimate in (\ref{local}) in the  Appendix we obtain
  \bqn
\frac{d}{dt} \|\bm \|_{L^p} &\leq& C(1+\|S\|_{\dot{B}^0_{\infty,
\infty}}) \log (1+\|S
\|_{W^{1,p}} ) \|\bm \|_{L^p} \quad (\mbox{for }\,\,p>N)\n \\
&\leq& C(1+\|S\|_{\dot{B}^0_{\infty, \infty}}) \log (1+\|D^2\bu
\|_{L^p} ) \|\bm \|_{L^p}\\
&\leq &C(1+\|S\|_{\dot{B}^0_{\infty, \infty}}) \log (1+\|\bm
\|_{L^p} ) \|\bm \|_{L^p}
 \eqn
 for $p>N$, where we used the boundedness on $L^p (\Bbb R^N)$ of the
 pseudo-differential operator
 $$\sigma_{ij} (D):=\partial_i\partial_j (1-\alpha \Delta)^{-1}
 = -R_iR_j \Delta (1-\alpha \Delta)^{-1}
 $$
 with the Riesz transforms $\{ R_j\}_{j=1}^N$ on $\Bbb R^N$(see
 Lemma 2.1, pp. 133\cite{ste}), which provides us with
  \bb\label{pseudo}
   \|D^2 \bu \|_{L^p}= \sum_{i,j=1}^N \| \sigma_{ij}(D) \bm \|_{L^p} \leq C
 \|\bm\|_{L^p}
  \ee
 for all $p\in (1, \infty)$.
  By Gronwall's lemma we obtain
  \bb\label{mest}
\log \left( 1+ \|\bm (t)\|_{L^p} \right)
\leq  \log (1+\|\bm_0\|_{L^p})
 \exp\left(C\int_0 ^t
 (1+ \|S(\tau)\|_{\dot{B}^0_{\infty, \infty}} )d\tau
 \right)
  \ee
  for $p>N$. This implies that
  \bb\label{BUC_p}
     \lim\sup_{t\to t_*} \|\bm (t)\|_{L^p} =\infty \quad \mbox{if and
  only if }\quad \int_0 ^{t_*} \|S (t)\|_{\dot{B}^0_{\infty, \infty}} dt =\infty \,.
  \ee

  {\it Step 2}. \quad
  Taking derivative of (EP) and  $L^2(\Bbb R^N )$ inner product it  with $D\bm|D\bm |^{p-2}$,
we find that
 \bqn
 \lefteqn{\frac{1}{p} \frac{d}{dt} \|D\bm (t)\|_{L^p}^p
 = \frac{1}{p} \int_{\Bbb R^N} (\div \bu ) |D \bm |^p \,
 dx-\int_{\Bbb R^N} (D\bu \cdot \nabla )\bm \cdot  D\bm |D\bm
 |^{p-2}\, dx}\\
 &&\quad - \int _{\Bbb R^N} D (\nabla \bu )^\top \bm \cdot D\bm |D
 \bm |^{p-2}\, dx - \int_{\Bbb R^N}(\nabla \bu )^\top D\bm \cdot D\bm
 |D\bm |^{p-2} \, dx\\
 &&\quad-\int_{\Bbb R^N}D (\div \bu )\bm \cdot D\bm
 |D\bm |^{p-2} \, dx -\int_{\Bbb R^N}(\div \bu)\, D\bm \cdot D\bm
 |D\bm |^{p-2} \, dx\\
&& \leq\left(3+\frac{1}{p} \right)\int_{\Bbb R^N} |D \bu||D\bm
|^p\,dx
+2\int_{\Bbb R^N}|D^2 \bu| |\bm ||D\bm |^{p-1}\, dx\\
&&\leq\left(3+\frac{1}{p} \right) \|D \bu\|_{L^\infty}\|D\bm
\|_{L^p}^p +2\|D^2
\bu\|_{L^{2p}} \|\bm \|_{L^{2p}}\|D\bm \|_{L^p} ^{p-1} \\
&&\leq C \|\bm\|_{L^p}\|D\bm \|_{L^p}^p +C \|\bm \|_{L^{2p}}^2\|D\bm
\|_{L^p} ^{p-1} \eqn
  for $p>N$,
  where we used the Sobolev embedding and (\ref{pseudo}) to estimate
  $$ \|D\bu \|_{L^\infty} \leq C \|D^2 \bu \|_{L^p}\leq
  C\|\bm\|_{L^p}
  $$
  for $p>N$. Hence, for $p>N$ we have
$$
  \frac{d}{dt} \|D\bm (t)\|_{L^p}\leq C\|\bm\|_{L^p}\|D\bm \|_{L^p} +C
  \|\bm
\|_{L^{2p}}^2.
$$
By Gronwall's lemma, we have
\bb \label{EDp}
 \|D\bm (t)\|_{L^p}\leq
 \exp\left( C\int_0 ^t  \|\bm (\tau)\|_{L^p} d\tau\right) \left(  \|D\bm_0 \|_{L^p} +C \int_0 ^t\|\bm(\tau ) \|_{L^{2p}}^2d\tau \right)
\ee
for $p>N$.
>From estimate (\ref{mest}), one has
%
\bq \label{e1p}
 & \int_0^t  \|\bm (s)\|_{L^p} ds
  \le   t   \max_{0\le s\le t} \|\bm (s)\|_{L^p}  \n \\
 &\quad
    \le   t   \max_{0\le s\le t}\exp \left(\log(1+ \|\bm (s)\|_{L^p}) \right) \n \\
&\quad  \leq  t \exp \left(   \log (1+\|\bm_0\|_{L^p})
 \exp\left(C\int_0 ^t
 (1+ \|S(\tau)\|_{\dot{B}^0_{\infty, \infty}} )d\tau
 \right)  \right).
\eq
Similarly,
\bb\label{e2p}
  \int_0^t  \|\bm (s)\|_{L^{2p} } ds
\leq   t \exp \left(   \log (1+\|\bm_0\|_{L^{2p}})
 \exp\left(C\int_0 ^t
 (1+ \|S(\tau)\|_{\dot{B}^0_{\infty, \infty}} )d\tau
 \right)  \right).
\ee
 Combining (\ref{EDp}, \ref{e1p}) and (\ref{e2p}), one obtains
%
  \bb\label{BUC_1p}
     \lim\sup_{t\to t_*} \|D \bm (t)\|_{L^p} =\infty \quad \mbox{if and
  only if }\quad \int_0 ^{t_*} \|S (t)\|_{\dot{B}^0_{\infty, \infty}} dt =\infty \,.
  \ee

{\it Step 3}. \quad
 Similarly, taking $D^2$ of (EP) and  $L^2(\Bbb R^N )$ inner
product it  with $D^2\bm|D^2\bm |^{p-2}$, we find that
 \bqn
 \lefteqn{\frac{1}{p} \frac{d}{dt} \|D^2\bm (t)\|_{L^p}^p
 \leq 4\int_{\Bbb R^N} |D \bu| |D^2 \bm |^p \, dx +
 3\int_{\Bbb R^N} |D^2 \bu ||D\bm| |D^2 \bm|^{p-1}\, dx}\\
 &&\qquad + 2\int_{\Bbb R^N} |D^3 \bu| |\bm||D^2 \bm |^{p-1}\, dx\\
&&\leq 4 \|D\bu\|_{L^\infty} \|D^2 \bm \|_{L^p}^p +3\|D^2
\bu\|_{L^{2p}} \|D \bm \|_{L^{2p}} \|D^2 \bm \|_{L^p}^{p-1} \\
&& \qquad+
2\|D^3 \bu\|_{L^{2p}} \|\bm\|_{L^{2p} } \|D^2 \bm \|_{L^{p}}^{p-1}\\
&&\leq C\|\bm \|_{L^p} \|D^2 \bm \|_{L^p}^p +C\|\bm\|_{L^{2p}}
\|D\bm \|_{L^{2p}} \|D^2 \bm \|_{L^p}^{p-1}
  \eqn
    for $p>N$, where we used the estimate (\ref{pseudo}) as follows
   \bqn \|D^3 \bu \|_{L^q}&=& \|\left\{ D^2 (1-\alpha \Delta )^{-1}\right\}
    D (1-\alpha \Delta ) \bu\|_{L^q}\\
     &\leq&  \sum_{i,j=1}^N
    \|\sigma_{ij}(D) D \bm \|_{L^q} \leq C \|D\bm \|_{L^q},
   \eqn
    which holds for all $q\in (1, \infty)$. Hence,
     $$
  \frac{d}{dt} \|D^2\bm (t)\|_{L^p}\leq C \|\bm \|_{L^p} \|D^2 \bm \|_{L^p}+C\|\bm\|_{L^{2p}}
\|D\bm \|_{L^{2p}}.
$$
 By Gronwall's lemma we have
$$
 \|D^2\bm (t)\|_{L^p}
 \leq
 \exp \left( C\int_0 ^t  \|\bm (\tau)\|_{L^p} d\tau\right)
 \left(   \|D^2 \bm _0 \|_{L^p} + C\int_0 ^t \|\bm(\tau)\|_{L^{2p}}
\|D\bm (\tau)\|_{L^{2p}}d\tau\right)
$$
for $p>N$. Similarly to the estimates in (\ref{e1p}) and (\ref{e2p}),
the right hand side terms in the above inequality
can all be controlled
$$
\int_0 ^t  (1+ \|S(\tau)\|_{\dot{B}^0_{\infty, \infty}} )d\tau \,.
$$
Therefore, we have
  \bb\label{BUC_2p}
     \lim\sup_{t\to t_*} \|D^2 \bm (t)\|_{L^p} =\infty \quad \mbox{if and
  only if }\quad \int_0 ^{t_*} \|S (t)\|_{\dot{B}^0_{\infty, \infty}} dt =\infty \,.
  \ee
Combination of (\ref{BUC_Hk}, \ref{BUC_p}, \ref{BUC_1p}, \ref{BUC_2p}) gives the proof of the theorem.
  $\square$

  \bigskip

  We now present a finite time blow-up result for $\alpha=0$.

 \begin{thm}
  Let the initial data of the system (\ref{hyp}), $\bu_0\in  H^k(\Bbb R^N)$, $k>N/2 +2$,  has the reflection
  symmetry with respect to the origin, and satisfies $\div \bu_0 (0) <0$. Then, there exists a
  finite time blow-up of the classical solution. 
 \end{thm}
{\it Proof } Taking divergence of (\ref{hyp}), we find
 \bb\label{diveq}
 \partial_t ( \div \bu)+\bu\cdot
\nabla(\div\bu) + 2\sum_{i,j=1}^N S_{ij}^2  +\sum_{j=1}^N (\Delta
u_j )u_j +(\div \bu)^2 + \sum_{i,j=1}^N (\partial_i\partial_j u_i)
u_j=0,
 \ee
  where we used $S_{ij} =\frac12 (\partial_i u_j +\partial_j
u_i ),$ and the fact
$$\sum_{i,j=1}^N
\partial_i u_j\partial_j u_i +\sum_{i,j=1}^N
\partial_i u_j\partial_i u_j= 2 \sum_{i,j=1}^N \partial_i u_j
S_{ij}=\sum_{i,j=1}^N (\partial_i u_j +\partial_j u_i ) S_{ij}=2
\sum_{i,j=1}^N  S_{ij}^2. $$
 Now we consider the reflection transform:
 $$ R: x\to \bar{x}=-x, \quad \bu(x,t)\to \bar{\bu}(x,t)=-\bu(-x,t).
 $$
 Obviously the system (\ref{hyp}) is invariant under this transform.
 The origin of the coordinate is the invariant point under the reflection
 transform. We consider the smooth initial data $\bu_0 \in H^k (\Bbb R^N)$, $k>N/2 +2$, which has the reflection symmetry. Then, by the
 uniqueness of the local classical solution in $H^k (\Bbb R^N)$, and hence in $C^2( \Bbb R^N  )$, the
 reflection symmetry is preserved as long as the classical solution
 persists.
 We consider the evolution of the solution at the origin of the
 coordinates.
 Then, $\bu(0,t)=0$ and  $D^2\bu(0,t)=0$ for all $t\in [0, T_*)$, where $T_*$ is the maximal time of existence of the classical
 solution in $H^k (\Bbb R^N)$.  If $T_* =\infty$, we will show that this leads to a
 contradiction.
  The system (\ref{diveq}) at the
 origin is reduced to
  $$
 \partial_t ( \div \bu)+ 2\sum_{i,j=1}^N S_{ij}^2  +(\div \bu)^2=0,
  $$
 which implies
  \bb
 \partial_t ( \div \bu)=-2\sum_{i,j=1}^N S_{ij}^2-(\div \bu)^2\leq -(\div \bu)^2.
 \ee
 Thus,
 $$
 \div \bu(0, t)\leq \frac{ \div \bu _0 (0)}{1+\div \bu_0 (0) t},
 $$
 which shows $T_* \leq \frac{1}{|\div \bu_0 (0)|}$ for $\div \bu_0 (0) <0$.
$\square$


\section{Zero  $\alpha$ limit for weak solutions}
In this section, we show the following theorem on the zero dispersion limit $\alpha \to 0$
for the weak solutions.
\begin{thm}
Let $\bu^\a \in  L^\infty( (0, T); H^1(\Bbb R^N))$ be a weak solution with initial data $\bu^\a_0$  to (EP) with $\a >0$,
and $\bu\in L^\infty( (0, T);H^k(\Bbb
R^N))\cap Lip ( (0, T); H^{2} (\Bbb R^N))$,  $k>N/2+3$, be the classic solution  with initial data $\bu_0$ to (EP)
with $\a=0$, i.e., (\ref{hyp}). Then, we have
\bq
&\sup_{0\leq t\leq T} \{ \|\bu^\a (t) -\bu(t)\|_{L^2}+ \sqrt{\a} \|\nabla (\bu^\a(t)- \bu (t)) \|_{L^2}\} \n  \\
& \qquad \leq
C \left( \a +
 \| \bu^\a_0 -\bu_0 \|_{L^2}
 + \sqrt{\a} \|\nabla ( \bu^\a_0-\bu_0 ) \|_{L^2}
\right),
\eq
where $C=C(\|u\|_{L^\infty(0,T;H^k(\Bbb R^N))}, \|u\|_{Lip(0,T; H^{2} (\Bbb R^N))})$ is a constant.
\end{thm}

{\it Proof }  We denote $\bm:= \bu-\a \Delta \bu$. Then (\bu, \bm) satisfy (EP) with a truncation term as below
 \bb\label{zero}
 \partial_t \bm +\div
(\bu \otimes \bm)+(\nabla \bu)^\top \bm =-\a \left\{ \Delta \bu_t
+\mathrm{div} (  \bu \otimes \Delta\bu) + (\nabla \bu)^\top \Delta
\bu \right\} \,.
 \ee
Subtracting  (\ref{zero})
from the first equation of (EP), and setting $ \bar \bm:=
\bm^\a-\bm$ and $\bar \bu:= \bu^\a -\bu, $ we find
 \bq\label{mbar}
&&\partial_t  \bar \bm +\mathrm{div} ( \bar  \bu \otimes  \bar
\bm)+\mathrm{div} ( \bar \bu \otimes \bm) +\mathrm{div} ( \bu
\otimes \bar \bm) +(\nabla \bar \bu)^\top \bar \bm +(\nabla \bar
\bu)^\top
\bm +(\nabla \bu)^\top \bar \bm\n\\
&&\qquad=\a \left\{ \Delta \bu_t +\mathrm{div} (\bu \otimes \Delta
\bu) + (\nabla \bu)^\top \Delta \bu \right\}
  \eq
Taking $L^2(\Bbb R^N )$ inner product (\ref{mbar}) with $\bar \bu$,
and integrating by part, and observing
$$
\int_{\Bbb R^N} \mathrm{div} ( \bar  \bu \otimes  \bar \bm)\cdot
\bar \bu \,dx=-\int_{\Bbb R^N}\bar \bu \cdot(\nabla \bar \bu)^\top
\bar \bm\, dx
$$
$$
\int_{\Bbb R^N}\mathrm{div} ( \bar \bu\otimes \bm)\cdot \bar \bu
\,dx=-\int_{\Bbb R^N}\bar \bu\cdot (\nabla \bar\bu)^\top  \bm \, dx,
$$
  we obtain that
 \bqn
 && \frac12 \frac{d}{dt}\int_{\Bbb R^N}\left(|\bar \bu|^2  +\a |\nabla \bar
 \bu|^2\right)\, dx=-\int_{\Bbb R^N}\mathrm{div} ( \bu \otimes \bar \bm)\cdot \bar \bu
\,dx-\int_{\Bbb R^N}\bar\bu\cdot (\nabla  \bu)^\top \bar\bm\, dx\\
&&\qquad +\a \int_{\Bbb R^N}\left[\bar \bu\cdot \left\{ \Delta
\bu_t +\mathrm{div} (  \bu \otimes \Delta\bu) + (\nabla \bu)^\top
\Delta \bu \right\}\right]\, dx\\
&&:=I_1+I_2+I_3. \eqn
  We estimate
  \bqn
  I_1&=& -\sum_{i,j=1}^N \int_{\Bbb R^N} \partial_i  u_i  (\bar u_j -\a \Delta \bar
  u_j )\bar u_j \, dx - \sum_{i,j=1}^N\int_{\Bbb R^N}  u_i\partial_i  (\bar u_j -\a \Delta \bar
 u_j )\bar u_j\, dx\n \\
&=&J_1+J_2, \eqn where
 \bqn
  J_1 &=&-\sum_{i,j=1}^N \int_{\Bbb R^N}\partial_i  u_i  |\bar u_j|^2
  dx+\a \sum_{i,j,k=1}^N \int_{\Bbb R^N} \partial_i\partial_k   u_i
  (\partial_k \bar u_j)\bar  u_j dx\\
  &&\qquad +\a
  \sum_{i,j,k=1}^N \int_{\Bbb R^N} \partial_i  u_i
  (\partial_k \bar u_j)\partial_k \bar  u_j dx\\
  &\leq & C \|u(t)\|_{C^2} (\|\bar \bu\|_{L^2}^2 +\a \|\nabla \bar \bu
  \|_{L^2}^2 ),
\eqn and
 \bqn
 J_2&=&- \sum_{i,j=1}^N\int_{\Bbb R^N} u_i(\partial_i  \bar u_j)\bar
 u_j\, dx +\a \sum_{i,j=1}^N\int_{\Bbb R^N}  u_i\partial_i (\Delta \bar
 u_j )\bar u_j \, dx\\
 &=&- \frac12 \sum_{i,j=1}^N\int_{\Bbb R^N}u_i \partial_i |\bar
 u_j|^2\, dx -\a \sum_{i,j,k=1}^N\int_{\Bbb R^N} \partial_k  u_i\partial_i (\partial_k \bar
 u_j )\bar u_j \, dx\\
 &&\qquad -\frac{\a}{2}\sum_{i,j,k=1}^N\int_{\Bbb R^N}  u_i\partial_i |\partial_k \bar
 u_j |^2 \, dx\\
 &=& \frac12\sum_{i,j=1}^N\int_{\Bbb R^N}\partial_i u_i  |\bar
 u_j|^2\, dx+\a \sum_{i,j,k=1}^N\int_{\Bbb R^N} \partial_i\partial_k  u_i (\partial_k \bar
 u_j )\bar u_j \, dx\\
 &&\qquad +\a \sum_{i,j,k=1}^N\int_{\Bbb R^N}\partial_k  u_i (\partial_k \bar
 u_j ) \partial_i\bar u_j \, dx
 +\frac{ \a}{2}\sum_{i,j,k=1}^N\int_{\Bbb R^N}  \partial_i u_i |\partial_k \bar
 u_j |^2 \, dx\\
 &\leq & C \|u(t)\|_{C^2} (\|\bar \bu\|_{L^2}^2 +\a  \|\nabla \bar \bu
  \|_{L^2}^2 ).
\eqn
 \bqn
 I_2&=&\sum_{i,j=1}^N\int_{\Bbb R^N}\bar u_i \partial_i u_j (\bar u_j
 -\a \Delta \bar u_j)
 \, dx\\
 &=&\sum_{i,j=1}^N\int_{\Bbb R^N}\bar u_i \partial_i u_j \bar
 u_j \, dx +\a \sum_{i,j,k=1}^N\int_{\Bbb R^N}\partial_k \bar u_i
 \partial_i u_j \partial_k \bar u_j \, dx\\
 &&\qquad +\a \sum_{i,j,k=1}^N\int_{\Bbb R^N} \bar u_i
 \partial_i \partial_k u_j \partial_k \bar u_j \, dx\\
 &\leq & C\|u(t)\|_{C^2} (\|\bar \bu\|_{L^2}^2 +\a  \|\nabla \bar \bu
  \|_{L^2}^2 ).
 \eqn
 One can estimate $I_3$ immediately as
  $$
 I_3 \leq \|\bar \bu \|_{L^2}^2 +\a^2C(\|\bu\|^2_{Lip(0, T;H^2(\Bbb
 R^N))} + \| \bu\|^4_{L^\infty (0, T; H^3(\Bbb R^N))} ).
 $$
 Summarizing the above estimates, we obtain
 $$
 \frac{d}{dt} (\|\bar \bu\|_{L^2}^2 +\a \|\nabla \bar \bu
  \|_{L^2}^2 )\leq C  \|u(t)\|_{C^2}  (\|\bar \bu\|_{L^2}^2 +\a \|\nabla \bar \bu
  \|_{L^2}^2 )+\a^2 C(\|\bu\|^2_{Lip(0, T;H^2(\Bbb
 R^N))} + \| \bu\|^4_{L^\infty (0, T; H^3(\Bbb R^N))} ),
  $$
which implies  by Gronwall's lemma  that
$$
 \|\bar \bu\|_{L^2}^2 +\a  \|\nabla \bar \bu \|_{L^2}^2
 \leq C_1( \a^2  +  \|\bar \bu(0) \|_{L^2}^2 +\a  \|\nabla \bar \bu(0) \|_{L^2}^2)
  $$
  where constant $C_1$ depended only on $\|\bu\|_{Lip(0, T;H^2(\Bbb
 R^N))}$ and $ \| \bu\|_{L^\infty (0, T; H^3(\Bbb R^N))}$.
  This completes the proof of theorem.
  $\square$

\section{ Liouville type theorem for stationary solutions}

In this section, we prove a Liouville type theorem for stationary solutions. Recall that the  stationary weak solution defined in Definition 1 reduces to

\begin{defin}
$\bu\in H^1(\Bbb R^N)$ is a stationary weak solution to
(EP) on $\Bbb R^N$, if the following holds
 \bq\label{defs}
&& \sum_{j=1}^N\int_{\Bbb R^N} \left\{u_i u_j+\alpha \nabla u_i
\cdot \nabla
 u_j \right\}\partial_j  \varphi_i \, dx +
 \a \sum_{j=1}^N\int_{\Bbb R^N}  u_j \nabla u_i\cdot \nabla \partial_j \varphi \n_i  \, dx \\
 &&\quad
+\sum_{j=1}^N\int_{\Bbb R^N} \left\{\frac{\delta_{ij} }{2} |\bu|^2
-\a
\partial_i \bu\cdot
\partial_j \bu +\frac{\a\delta_{ij}}{2} |\nabla \bu|^2\right\}
\partial_j \varphi _i\, dx=0
\eq
 for $i=1, \cdots, N$ and for all $\bphi\in C_0 ^\infty (\Bbb R^N)$.
\end{defin}
\begin{thm}
(i) Let $\bu\in H^1(\Bbb R^N)$ be a stationary weak solution to (EP)
with $\alpha > 0$. Then, $\bu=0$.

(ii) Let $\bu\in L^{2}(\Bbb R^N)$ be a stationary weak solution to (EP)
with $\alpha = 0$. Then, $\bu=0$.
\end{thm}

\noindent{\it Proof } For $\alpha>0$, one can write (\ref{defs}) in the following
form,
 \bb\label{def2}
\sum_{j=1}^N\int_{\Bbb R^N} T_{ij}^a\partial_j \varphi _i\,
dx+\sum_{j,k=1}^N\int_{\Bbb R^N} \tilde{T}_{ijk}^b\partial_j\partial_k
\varphi_i \, dx=0, \ee
where $T_{ij}^a$ is defined in (\ref{Ta}) and we recall here
$$
T_{ij}^a=u_i u_j+\alpha \nabla u_i \cdot \nabla
 u_j +\frac{\delta_{ij} }{2} |\bu|^2
-\a
\partial_i \bu\cdot
\partial_j \bu +\frac{\a\delta_{ij}}{2} |\nabla \bu|^2,
$$
and
$$
\tilde{T}_{ijk}^b= \alpha u_j \partial_k u_i.
$$
corresponding to $T_{ij}^b$ in (\ref{Tb}).

 Let us consider the radial cut-off function
$\sigma\in C_0 ^\infty(\Bbb R^N)$ such that
 $$
   \sigma(|x|)=\left\{ \aligned
                  &1 \quad\mbox{if $|x|<1$},\\
                     &0 \quad\mbox{if $|x|>2$},
                      \endaligned \right.
$$
and $0\leq \sigma  (x)\leq 1$ for $1<|x|<2$.  Then, for each $R
>0$, we define
 $$
\s \left(\frac{|x|}{R}\right):=\s_R (|x|)\in C_0 ^\infty (\Bbb R^N).
$$
Choosing  $\varphi_i (x)=x_i \si$ in (\ref{def2}),  we obtain
 \bq\label{l2}
 0&=&\sum_{i=1}^N \int_{\Bbb R^N} T_{ii}^a\si \,dx+\sum_{i,j=1}^N
 \int_{\Bbb R^N} T_{ij}^a x_j \partial_i \si \, dx
 +\sum_{i,k=1}^N
 \int_{\Bbb R^N} \tilde{T}_{iik}^b\partial_k \si \, dx  \n \\
 && +\sum_{i,j=1}^N
 \int_{\Bbb R^N} \tilde{T}_{iji}^b\partial_j \si \, dx
  +\sum_{i,j,k=1}^N
 \int_{\Bbb R^N} \tilde{T}_{ijk}^b x_i \partial_j\partial_k \si \, dx \n \\
 &=&I_1+I_2+I_3 +I_4 + I_5.
 \eq
 The hypothesis $u\in H^1 (\Bbb R^N)$ implies that $T \in L^1 (\Bbb
 R^N)$.
Thus, we obtain
 $$
 |I_2|\leq \frac{1}{R} \int_{\{ R\leq |x|\leq 2R\} } |T^a| |x|
 |\nabla \s|\, d\bx
 \leq 2\|\nabla \s \|_{L^\infty} \int_{\{R\leq |x|\leq 2R\} } |T^a|\,
d\bx\to 0
 $$
 as $R\to \infty$ by the dominated convergence theorem.
Similarly, $I_3, I_4, I_5 \to 0$ as $R\to \infty$.

 Thus, passing
 $R\to \infty$ in (\ref{l2}), we have
 \bqn
0&=&\lim_{R\to \infty} \sum_{i=1}^{N} \int_{\Bbb R^N} T_{ii}^a \si\, d\bx \n \\
 &=&\int_{\Bbb R^N} \left\{ \frac{(N+2)}{2} |\bu|^2 +\frac{\a
 N}{2} |\nabla \bu|^2\right\} d\bx,
 \eqn
 which implies $\bu=0$. This gives (i).

For the case $\alpha=0$. All the terms involving $\alpha$ drop and (ii) holds true. This completes the proof the theorem
$\square$

We remark that the Liouville
type results in Theorem 5 is rather surprising, as all the previous
 Liouville
type results are for dissipative systems. For  the Liouville type results for the  dissipative systems, see, e.g. \cite{Chae}.
Theorem 5  is the first  Liouville
type theorem for non-dissipative systems.

\bigskip

\noindent {\bf Acknowledgements:} This work was initiated  at Duke
University when the first author visited there. The authors wish to
acknowledge the hospitality of  Mathematical Sciences Center of
Tsinghua University  where this research was completed. The research
of D.C was supported partially by  NRF Grant no. 2006-0093854. The
research of J.-G. L. was partially supported by NSF grant DMS
10-11738.

\section*{Appendix: proof of Theorem 1}

\noindent{\it Proof of Theorem 1  } The proof of local existence
part is standard, and below we derive the key local in time estimate
of $\bu (t) \in L^\infty ([0, T); W^{2,p} (\Bbb R^N))\cap Lip(0, T;
W^{1,p} (\Bbb R^N))$.
 \bq\label{local}
 \frac{1}{p}\frac{d}{dt} \|\bm \|_{L^p} ^p &=&\frac{1}{p}\int_{\Bbb R^N} (\bu
 \cdot \nabla ) |\bm |^p \,dx +\sum_{i,j=1}^N
 \int_{\Bbb R^N} \partial_j u_i m_i m_j |\bm |^{p-2} \, dx\n\\
 &&\qquad +\int_{\Bbb R^N } (\div \bu )|\bm |^p \, dx\n\\
 &=&\left(1-\frac{1}{p} \right) \int_{\Bbb R^N} Tr(S)|\bm |^p \, dx
 +\sum_{i,j=1}^N \int_{\Bbb R^N} S_{ij} m_im_j |\bm |^{p-2} \, dx\n\\
 &\leq & C
 \|S\|_{L^\infty} \|\bm \|_{L^p} ^p\leq  C
 \|\nabla \bu\|_{L^\infty} \|\bm \|_{L^p} ^p\leq  C
 \|\bm \|_{L^p} ^{p+1},
 \eq
 and therefore
 $$  \frac{d}{dt} \|\bm \|_{L^p}\leq C \|\bm \|_{L^p}^2. $$
We thus have the following estimate on $L^\infty (0, T; W^{2,p}
(\Bbb R^N))$,
 \bb
 \|\bu (t)\|_{W^{2,p}} \leq \frac{C\|\bu_0 \|_{W^{2,p}}}{1-C t\|\bu_0
 \|_{W^{2,p}}} \quad \forall t\in [0, T),
\ee
 where $T= \frac{1}{\|\bu_0 \|_{W^{2,p}}}$.
In order to have estimate of $\bu$ in $Lip(0, T; W^{1,p} (\Bbb
R^N))$, we take $L^2(\Bbb R^N)$ inner product (EP) with the test
function $\bpsi \in W^{1,\frac{p}{p-1}} (\Bbb R^N)$ for $p>N$. Then,
 \bqn
 \int_{\Bbb R^N }\partial_t \bm \cdot \bpsi \, dx
 &=& \int_{\Bbb R^N} \bm (\bu \cdot \nabla )\bpsi \, dx -\int_{\Bbb R^N}
 \bm\cdot \nabla \bu \cdot \bpsi\, dx   
 \\
 &\leq & C\|\bm \|_{L^p} \|\bu\|_{L^\infty} \|\nabla \bpsi\|_{L^{\frac{p}{p-1}}}+C\|\bm \|_{L^p} \|\nabla \bu\|_{L^\infty} \|\nabla
 \bpsi\|_{L^{\frac{p}{p-1}}}\n \\
 &\leq & C \|\bm \|_{L^p}^2 \|\bpsi\|_{W^{1,\frac{p}{p-1}}},
 \eqn
 which provides us with the estimate,
 $$
  \|\partial_t \bu\|_{L^\infty (0, T; W^{1,p}(\Bbb R^N))}
 \leq C \|\partial_t \bm\|_{L^\infty (0, T; W^{-1,p}(\Bbb R^N))}
 \leq C \|\bm \|_{L^\infty(0,T);L^p (\Bbb R^N)}^2.
 $$
 Hence, for all $0<t_1<t_2 <T$ we have
 $$
 \|\bu(t_2)-\bu(t_1)\|_{W^{1,p}}\leq \int _{t_1}^{t_2}\left\|\partial_t \bu
 (t)\right\|_{W^{1,p}} dt\leq  C(t_2-t_1)\|\bm \|_{L^\infty(0,
 T;L^p (\Bbb R^N))}^2.
 $$
Namely,
$$ \|\bu\|_{Lip(0, T; W^{1,p} (\Bbb
R^N))}\leq C\|\bm \|_{L^\infty(0,
 T;L^p (\Bbb R^N))}^2.
 $$
This gives (i). Next we prove local in time persistency of regularity for  $\bu(t)$
in $ H^{k} (\Bbb R^N))$ with $k>N/2+3$. Let $\beta=(\beta_1, \cdots,
\beta_N)$ be the standard multi-index notation with $|\beta|=\beta_1
+\cdots +\beta_N$. Taking $H^k(\Bbb R^N )$ inner product (EP) with
$\bm$, we find
 \bq \label{I123}
  \frac12 \frac{d}{dt}\sum_{|\beta|\leq k } \|D^\beta \bm\|_{L^2}^2 &=&
  -\sum_{|\beta|\leq k } \int_{\Bbb R^N} D^\beta\{ (\bu\cdot \nabla ) \bm\}\cdot D^\beta \bm
  \, dx  \n \\
  &&-\sum_{|\beta|\leq k } \int_{\Bbb R^N} D^\beta\{ (\nabla
\bu)^\top \bm\}\cdot D^\beta \bm
  \, dx  \n \\
  &&-\sum_{|\beta|\leq k } \int_{\Bbb R^N} D^\beta\{ (\mathrm{div}\,
 \bu)\bm\}\cdot D^\beta \bm
  \, dx \n \\
  &&:=I_1+I_2+I_3.
 \eq
We write
 \bqn
 I_1&=&-\sum_{|\beta|\leq k } \int_{\Bbb R^N} \{D^\beta(\bu\cdot \nabla )
  \bm- (\bu \cdot\nabla ) D^\beta\bm\}\cdot D^\beta \bm
  \, dx\\
  &&\qquad+\sum_{|\beta|\leq k } \int_{\Bbb R^N}(\bu \cdot\nabla ) D^\beta \bm \cdot D^\beta \bm\,
  dx\\
  &&:=J_1+J_2,
  \eqn
 and using the standard commutator estimate, we deduce
 \bq \label{J1}
 J_1&\leq&\sum_{|\beta|\leq k } \|D^\beta(\bu\cdot \nabla )
 \bm- (\bu \cdot\nabla ) D^\beta \bm\|_{L^2} \|D^\beta \bm \|_{L^2} \n \\
 &\leq &C(\|\nabla\bu \|_{L^\infty} \|\bm \|_{H^k} +\|\bu\|_{H^k} \|\nabla
 \bm\|_{L^\infty} )\|\bm\|_{H^k}\\
 &\leq & C (\|\bu\|_{H^{N/2 +1+\vare}} \|\bm\|_{H^k}
 +\|\bm\|_{H^{k-2}}
 \|\bm\|_{H^{N/2 +1+\vare}})\|\bm\|_{H^k}\, \,\,(\forall \vare >0) \n \\
 &\leq & C\|\bm\|_{H^k}^3 \n
 \eq
 if $k>N/2+1$, where we used the fact $\bu=(1-\a  \Delta )^{-1}\bm$, and therefore
$ \|\bu\|_{H^s}\leq \|\bm\|_{H^{s-2}}$ for all $s\in \Bbb R$.
 \bq\label{J2}
J_2& =&\frac12\sum_{|\s|\leq k } \int_{\Bbb R^N}(\bu \cdot\nabla )
|D^\s \bm|^2\,
  dx=-\frac12\sum_{|\s|\leq k } \int_{\Bbb R^N}(\div \bu ) |D^\beta
\bm|^2\,
  dx \n \\
  &\leq& C\|\nabla \bu\|_{L^\infty} \|\bm \|_{H^k}^2 \leq
  C\|\bm\|_{H^{N/2-1+\vare}}\|\bm \|_{H^k}^2\, \,\,(\forall \vare >0)  \\
&\leq& C \|\bm\|_{H^k}^3 \n
  \eq
if $k>N/2-1$. The estimates of $I_2,I_3$ are simpler, and we have
\bq \label{I23}
  I_2+I_3 &\leq &\|(\nabla \bu)^\top \bm\|_{H^k} \|\bm \|_{H^k}\leq
C(\|\nabla \bu\|_{L^\infty} \|\bm\|_{H^k} +\|\bu\|_{H^{k+1}}\|\bm
\|_{L^\infty} )\|\bm\|_{H^k}\\
&\leq & C (\|\bm \|_{H^{N/2-1+\vare}}\|\bm\|_{H^k}+
\|\bm\|_{H^{k-1}} \|\bm \|_{H^{N/2+\vare}} )\|\bm\|_{H^k }  \n \\
&\leq & C \|\bm \|_{H^k} ^3 \n
\eq
 if $k>N/2$. Summarizing the above estimates,  we obtain
 $$
 \frac{d}{dt} \|\bm \|_{H^k}^2 \leq C \|\bm \|_{H^k} ^3
 $$
 for $k>N/2+1$,
 which  implies
 $$ \|\bu (t) \|_{H^k} \leq \frac{C\|\bu_0\|_{H^k}}{1-C \|\bu_0
 \|_{H^k} t} \qquad \forall t\in [0, T),\,\mbox{where}\,\,
 T=\frac{1}{C\|\bu_0\|_{H^k}},
 $$
where $k>N/2+3$. \\
We now prove uniqueness of solution in this
class. Let $(\bu_1, \bm_1), (\bu_2, \bm_2) $ two solution pairs
corresponding to initial data $(\bu_{1,0}, \bm_{1,0}), (\bu_{2,0},
\bm_{2,0})$. We set $\bu=\bu_1-\bu_2$, and so on. Subtracting the
equation for $(\bu_2, \bm_2) $ from that of  $(\bu_1, \bm_1)$, we
find that
 \bb\label{dif}
 \partial_t \bm + \div \big(\bu_1 \otimes \bm\big) +\div \big(\bu \otimes \bm_2\big) +
 (\nabla \bu_1)^\top\bm + (\nabla \bu )^\top\bm_2 =0.
 \ee
 Let $p>N$. Taking $L^2 (\Bbb R^N)$ product of (\ref{dif}) with $\bm
 |\bm|^{p-2}$, we obtain
 \bqn
 \lefteqn{\frac{1}{p}\frac{d}{dt} \|\bm(t)\|_{L^p}^p=
 -\left(1-\frac{1}{p}\right) \int_{\Bbb R^N} (\div \bu_1 )|\bm|^p \, dx
 -\int_{\Bbb R^N} (\div \bu)\bm_2 \cdot \bm |\bm|^{p-2}\, dx}\n \\
 &&
-\int_{\Bbb R^N} (\bu\cdot \nabla )\bm_2\cdot \bm |\bm|^{p-2}\, dx
-\int_{\Bbb R^N}(\nabla \bu_1)^\top\bm \cdot \bm |\bm|^{p-2}\, dx\n \\
&&\qquad-\int_{\Bbb R^N}(\nabla \bu )^\top\bm_2\cdot \bm |\bm|^{p-2}\, dx\n \\
&&\leq C(\|\div \bu_1 \|_{L^\infty} \|\bm\|_{L^p}^p +\|\nabla
\bu\|_{L^\infty} \| \bm_2\|_{L^p}\|\bm\|_{L^p}^{p-1} + \|\bu\|_{L^p}
\|\nabla \bm_2 \|_{L^\infty} \|\bm\|_{L^p}^{p-1}\n \\
&& \qquad + \|\nabla \bu_1 \|_{L^\infty} \|\bm\|_{L^p} ^{p} +
\|\nabla \bu\|_{L^\infty} \|\bm_2 \|_{L^p} \|\bm\|_{L^p} ^{p-1} )\n
\\
&& \leq C (\|\bu_1\|_{H^k} +\|\bu_2 \|_{H^k} ) \|\bm\|_{L^p} ^p
 \eqn
for $k>N/2+3$. Hence,
$$
\|\bm(t)\|_{L^p}\leq \|\bm_0\|_{L^p}\exp\left(C \int_0 ^t(\|\bu_1
(\tau)\|_{H^k} +\|\bu_2 (\tau)\|_{H^k} ) d\tau\right).
$$
This inequality implies the desired uniqueness of solutions in the
class $ L^1(0, T; H^k (\Bbb R^N))$ with $k>N/2+3$. This gives (ii).
The proof of (iii) was explained at the end of Section 2. This completes the proof of Theorem 1.
   $\square$

\bigskip


 \begin{thebibliography}{1}

 \bibitem{Arnold} Arnold, V.,
{\it Sur un principe variationnel pour les ecoulements stationnaires des liq- uides parfaits et ses applications aux probl`emes de stanbilit?e non lin?eaires},
J. M\'ec., {\bf 5} (1966),  29--43.

 \bibitem{Bardos} Bardos, C.,  Linshiz, J., and Titi, E.S.,
 {\it  Global regularity and convergence of a Birkhoff-Rott-$\alpha$ approximation of the dynamics of vortex sheets of the 2D Euler equations},
  Comm.  Pure and Appl. Math.,
 {\bf  63}  (2010), 697--746.

 \bibitem{BrCo2007} Bressan, A. and Constantin, A.,
{\it Global conservative solutions of the {C}amassa-{H}olm equation},
Arch. Ration. Mech. Anal., {\bf 183} (2007), 215--239,

 \bibitem{CH} Camassa, R. and Holm, D.D.,
{\it An integrable shallow water equation with peaked solitons},
Phys. Rev. Lett., {\bf 71} (1993), 1661--1664.

 \bibitem{Chae} Chae, D.,
 {\it On the nonexistence of global weak solutions to the Navier-Stokes-Poisson equations in $\Bbb R^N$},
 Comm. PDE, {\bf 35}  (2010), 535--557.

  \bibitem{Chertock} Chertock, A., Liu, J.-G., and Pendleton, T.,
{\it  Convergence of a particle method and global weak solutions
for a family of evolutionary PDEs},  submitted.

\bibitem{Chen} Chen, S., Foias, C., Holm, D.D., Olson, E., Titi, E.S., and Wynne, S.,
{\it Camassa-{H}olm equations as a closure model for turbulent channel and pipe flow},
Phys. Rev. Lett., {\bf 81} (1998), 5338--5341.

 \bibitem{CE98a} Constantin, A. and Escher, J.,
 {\it Global weak solutions for a shallow water equation},
Indiana Univ. Math. J., {\bf 47} (1998), 1527--1545.

\bibitem{CE98c} Constantin, A. and Escher, J.,
 {\it Global existence and blow-up for a shallow water equation},
Ann. Scuola Norm. Sup. Pisa Cl. Sci.  Serie IV, {\bf 26} (1998), 303--328.


\bibitem{CM} Constantin, A. and Molinet, L.,
 {\it Global weak solutions for a shallow water equation},
Comm. Math. Phys., {\bf 211} (2000), 45-61.

 \bibitem{EbinMarsden}
Ebin, D.  and Marsden, J.,
{\it Groups of diffeomorphisms and the motion of an incompressible fluid},
Ann. of Math, {\bf 92} (1970), 102--163.

\bibitem{Foias} Foias, C., Holm D.D., and Titi, E.S.,
{\it The three dimensional
viscous Camassa-Holm equations, and their relation to the Navier-Stokes equations and turbulence theory},
J. Dyn. and Diff. Eqns., {\bf 14} (2002), 1--35.

\bibitem{Fokas} Fuchssteiner, B. and Fokas, A.S,
 {\it Symplectic structures, their {B}\"{a}cklund transformations and hereditary symmetries},
Physica D, {\bf 4} (1981), 47--66.

\bibitem{HiMaAr2001} Hirani, A.N., Marsden, J.E. and Arvo, J.,
 {\it Averaged template matching equations},
Lecture Notes in Computer Science, volume 2134, EMMCVPR,
Springer,
(2001), 528--543.


\bibitem{Khesin} Khesin, B and Wendt, R,
The geometry of infinite-dimensional groups,
Springer, 2009.

\bibitem{Majda} Majda, A.,
{\it Compressible Fluid Flow and Systems of Conservation Laws in Several Space Variables},
Springer-Verlag 1984.

\bibitem{HoMa2005} Holm, D.D. and Marsden, J.E.,
{\it Momentum maps and measure-valued solutions (peakons,
              filaments, and sheets) for the {EPD}iff equation},
The breadth of symplectic and Poisson geometry,
Progr. Math., VOLUME 232, Birkh\"auser Boston,  (2005) 203--235.

 \bibitem{HMR1} Holm, D.D.,  Marsden, J.E. and Ratiu, T.S.,
{\it Euler-Poincar\'{e} models of ideal fluids with nonlinear dispersion},  Phys. Rev. Lett., {\bf 80} (1998), 4173--4177.


\bibitem{HMR2} Holm, D.D., Marsden, J.E. and Ratiu, T.S.,
{\it Euler-Poincar\'{e}  equations and semi-direct products
with applications to continuum theories},
Adv. in Math., {\bf 137} (1998),  1--81.

\bibitem{HNP} Holm, D.D., Nitsche, M., and Putkaradze, V.,
{\it Euler-alpha and vortex blob regularization of
vortex filament and vortex sheet motion},
J. Fluid Mech., {\bf 555}  (2006), 149--176.

\bibitem{HoRaTrYo2004} Holm, D.D., Ratnanather, J.T., Trouv\'{e}, A. and Younes, L.,
{\it Soliton dynamics in computational anatomy},
NeuroImage, {\bf 23} (2004), S170 - S178.


\bibitem{HSS} Holm, D.D., Schmah, T. and  Stoica, C.,
Geometric Mechanics and Symmetry: From Finite to Infinite Dimensions , Oxford University Press, 2009.

\bibitem{TitiXin}Jiu, Q.S., Niu, D.J. ,  Titi, E.S. and Xin, Z.P.,
{\it The Euler-$\alpha$ approximations to the 3D axisymmetric Euler
equations with vortex-sheets initial data}, preprint, (2009).

\bibitem{McKean} McKean, H.P.,
{\it Breakdown of the Camassa-Holm equation},
Comm. Pure Appl. Math., {\bf 57} (2004),  416--418.


\bibitem{MoLu2004} Molinet, L.,
 {\it On well-posedness results for {C}amassa-{H}olm equation on the
              line: a survey},
  J. Nonlinear Math. Phys., {\bf 11} (2004), 521--533.

\bibitem{Marsden} Marsden, J.E. and Shkoller, S.,
{\it Global well-posedness for the
Lagrangian averaged Navier-Stokes (LANS-$\alpha$) equations on bounded
domains},
Proc. Roy. Soc. London A, {\bf 359} (2001), 1449--1468.

\bibitem{ste} Stein, E.M.,
{\it Singular Integrals and Differentiability
Properties of Functions,}, Princeton, NJ,  Princeton
University Press. 1970.

\bibitem{tay} Taylor, M.,
{\it Tools for PDE}, AMS Mathematical
Surveys and Monographs {\bf 81}, (2000).

 \bibitem{XZ01} Xin, Z. and Zhang, P.,
  {\it On the weak solutions to a shallow water equation},
 Comm. Pure Appl. Math., {\bf 53} (2000), 1411--1433.

 \bibitem{Younes} Younes, L.,
Shapes and Diffeomorphisms,
Springer,  2010.


\end{thebibliography}
\end{document}